\documentclass[a4paper,11pt,leqno]{amsart}
\usepackage{amsmath,amscd,amsfonts,amssymb}
\usepackage{latexsym}
\usepackage{graphics}
\usepackage{epsfig}
\usepackage{tikz}
\usepackage[english]{babel}
\usetikzlibrary{calc}

\usepackage{ragged2e}
\usepackage{fancyhdr}
\numberwithin{equation}{section}
\newtheorem{theorem}{Theorem}[section]
\newtheorem{prop}[theorem]{Proposition}

\theoremstyle{definition}
\newtheorem{defi}{Definition}

\newtheorem{rem}{Remark}
\newtheorem{conj}{Conjecture}

\newcommand{\RR}{\mathbb{R}}
\newcommand{\ZZ}{\mathbb{Z}}
\newcommand{\KK}{\mathbb{K}}
\newcommand{\PP}{\mathbb{P}}

\title{Bavard's systolically extremal Klein bottles and three dimensional applications}

\author{Chady El Mir}
\email{chady.mir@gmail.com}
\address{Current address : \newline
Universit\'e Libanaise, Laboratoire de Math\'ematiques et Applications (LaMA), Tripoli, Liban}

\graphicspath{{images/}}

\begin{document}

\begin{flushright}

\textit{"In the memory of a dedicated mathematician and dear friend,}\\
\textit{Jacques Lafontaine."}

\end{flushright}

\date{}

\selectlanguage{english}
\begin{abstract}
A compact manifold is called Bieberbach if it carries a flat Riemannian metric. Bieberbach manifolds satisfy an isosystolic inequality by a general and fundamental result of M. Gromov. In dimension $3$, there exist four classes of non-orientable Bieberbach manifolds up to an affine diffeomorphism. In this paper, We prove the existence on each diffeomorphism class of non-orientable Bieberbach $3$-manifolds of a two-parameter family of singular Riemannian metrics that are systolically extremal  in their conformal class. The proof uses a one-parameter family of singular Riemannian metrics on the Klein bottle discovered by C. Bavard (\cite{bavard88}): each one of these metrics is extremal in its conformal class.   
\end{abstract}

\medskip

\maketitle

\selectlanguage{english} \noindent {\bf Keywords} : Systole, systolic
volume, singular Riemannian metric, Bieberbach
manifold, Radon measure,.
\\
\newline

\section{Introduction and main results}\label{intro}
The systole of a compact connected but non simply connected Riemannian manifold $(M,g)$ of dimension $n$ is the least length of a non-contractible curve and a systolic geodesic is a closed non trivial curve which realizes the systole. An isosystolic inequality on the manifold $M$ is an inequality of the form 
$$
\frac{\mathrm{Sys}(g)^n}{\mathrm{Vol}(g)}\leq c < +\infty$$

\noindent which holds for any Riemannian metric $g$ on $M$. This inequality leads us to introduce the systolic volume (or area in the two dimensional case) of $M$ defined as

$$\sigma(M):=
\sup_g \frac{\mathrm{Sys}(g)^n}{\mathrm{Vol}(g)}$$

\noindent where $g$ runs over all the Riemannian metrics on $M$.

The first result concerning isosystolic inequalities and marking the birth of the systolic geometry is due to C. Lowner. He showed that $\sigma(T^2)=\frac{2}{\sqrt{3}}$ with the supremum being attained by the flat equilateral torus (up to homothety). Following this direction, his student P. M. Pu, showed (cf. \cite{pu}) that $\sigma(\RR \mathbb{P}^2)=\frac{\pi}{2}$ with the supremum being attained by the round metric. In 1986, C. Bavard studied the case of the Klein bottle (cf. \cite{bavard86}) and showed that $\sigma(T^2)=\frac{2}{\sqrt{3}}$. This time the extremal metric, i.e. the metric that achieves the supremum, is not of constant curvature: it is of constant curvature equal to +1 outside a singular line (see Section \ref{klein}). The remaining surfaces of genus $g\geq 2$ also satisfy an isosystolic inequality but a sharp bound is still unknown for any genus (see \cite{hebda} and \cite{calabi}).
For more details about the systolic geometry, the reader is referred to the references \cite{bki} and \cite{katz} (also check \cite{katsab} and the references therein).\\


In dimensions $\geq 3$, a deep result due to M. Gromov asserts that every essential manifold (in particular aspherical manifolds) satisfies an isosystolic inequality (see \cite{gromov} for a definition of essential manifolds). Nevertheless, apart from this result, very little things are known about the systolic volume and the metrics that
can realize it. For example, it is not known in the simple cases of  tori and real projective spaces, whether the metrics of constant curvature are extremal. In a paper of 2008 (\cite{elaf08}), J. Lafontaine and the author studied the systolic geometry Bieberbach manifolds, i.e., compact connected manifolds that carry a flat Riemannian metric. Bieberbach manifolds are essential and satisfy then an isosystolic inequality. They showed that flat Riemannian metrics cannot be extremal on non-orientable Bieberbach 3-manifolds. It is then justified to search for extremal metrics on these manifolds among those that are not flat, preferably with ``many" systolic geodesics (cf. \cite{calabi}, also check Section \ref{maximus}). In the present paper, we show that there exists on each non-orientable Bieberbach 3-manifold a two parameter family of singular Riemannian metrics each of which is extremal in its conformal class. As a corollary, we deduce the result of \cite{elaf08}.\\

\begin{theorem} There exists on each non-orientable Bieberbach $3$-manifolds a two-parameter family of singular (non constant curvature) Riemannian metrics each of which is extremal in its conformal class.

\end{theorem}

The ingredients of our proof come from a paper of C. Bavard entitled ``In\'egalit\'e isosystolique conforme pour la bouteille de Klein" (\cite{bavard88}). In that paper of 1988, Bavard exhibits for each conformal class of the Klein bottle a Riemannian metric which is conformally extremal i.e. which is extremal in its conformal class. The resulting metrics are singular and locally either flat or spherical (see Section \ref{klein}). His proof (and ours) relies on the following result (initiated by Jenkins \cite{jenk}).\\

Let $(M,g)$ be a compact essential Riemannian manifold of dimension $n$ and let $\Gamma$ be the space of the systolic geodesics of $(M,g)$. For every Radon measure $\mu$ on $\Gamma$, we associate a measure $^* \mu $ on $M$ by setting for $\varphi \in C^0(M,\RR)$\\

$$ <^* \mu, \varphi>= <\mu,\overline{\varphi}>$$

\noindent where $\overline{\varphi}(\gamma)=\int{\varphi \circ \gamma (s) ds}$, $ds$ is the arc length of $\gamma$ with respect to $g$. Then we have

\begin{theorem}\label{confext} (\cite{bavard88},\cite{bavard92} and \cite{gromov}) The Riemannian manifold $(M,g)$ is conformally extremal if and only if there exists a positive measure $\mu$, of mass $1$, on $\Gamma$ such that

$$^*\mu=\frac{Sys(g)}{Vol(g)} \cdot dg$$

where $dg$ is the volume measure of $(M,g)$.

\end{theorem}

\begin{rem} The previous theorem can be generalized to compact length spaces with a general notion of volume (cf. \cite{bavard92}).

\end{rem}


\emph{Acknowledgements:} The author thanks the Lebanese Association for Scientific Research for hospitality and
support. It is a pleasure to thank Karim Youssef for his help in the realization of the figures and the use of MATLAB.

\section{Isosystolic inequalities on the Klein bottle}\label{klein}

A flat Klein bottle is the quotient  $\mathbb{R}^2/\Gamma$, where
$\Gamma$ is the subgroup of isometries of $\mathbb{R}^2$ generated by the glide reflection  $(x,y)\mapsto
(x+\frac{a}{2},-y)$ and the translation $(x,y)\mapsto (x,y+b)$. The connected component of the identity in $\mathrm{Isom}(\KK)$ consists of the horizontal translations
$$r_\delta :(x,y)\mapsto (x+\delta,y)$$
where $\delta$ is modulo $a$. The quotient $\mathrm{Isom}(\KK)/\mathrm{Isom}_0(\KK)$ is isomorphic to the \emph{Klein group}. The three non-trivial elements of this quotient can be represented by

\begin{enumerate}
\item the translation (after going to the quotient) $(x,y)\mapsto (x,y+ b/2)$. We denote by $T$ this transformation.

\item a reflection with respect to a vertical geodesic, which is also a symmetry with respect to a point on a horizontal short geodesic. We denote by $S_1$ such a transformation.

\item a symmetry with respect to a point on the common boundary of the two M\"obius bands composing the Klein bottle. We denote by $S_2$ such a symmetry.

\end{enumerate}

We will often consider flat Klein bottles that are the quotient of
$ \mathbb{C}$ by
the group generated by the maps
$\tau: z\mapsto \bar{z}+\pi$ and $t_{v}:z\mapsto z+2i\beta$. The induced flat metric on the Klein bottle will be denoted by $ g_\beta$. By the uniformization theorem, every Riemannian metric $g$ on the Klein bottle is conformally equivalent to a unique flat metric $ g_\beta$, for some $\beta \in (0,+\infty)$. The parameter $\beta$ represents then the conformal type of the metric $g$.\\


\noindent We denote by $C^{0,\infty}(\mathbb{R})$ the set of continuous and piecewise smooth functions on $\mathbb{R}$. \\

Suppose now that $h^\beta$ is a Riemannian metric on $\mathbb{C}\simeq \mathbb{R}^2$ which is periodic with respect to the fundamental domain $[-\frac{\pi}{2},\frac{\pi}{2}]\times[-\beta,\beta]$ and satisfying $(h^\beta)_{(x,y)}=\varphi(y)(dx^2+dy^2)$, where $\varphi\in C^{0,\infty}(\mathbb{R})$ is a positive, even and $2\beta$-periodic function. Then $h^\beta$ can be written in the form $(h_b)_{(u,v)}=f^2(v)du^2+dv^2$, where  $f\in C^{0,\infty}(\mathbb{R})$ is positive, even and $2b$-periodic function. To see this, let $\phi(y)=\int_0^{y}{\sqrt{\varphi(t)}dt}$ and let $f=(\varphi \circ \phi ^{-1})^{\frac{1}{2}}$. The map $G:(x,y) \rightarrow (x,\phi(y))$ is then a diffeomorphism on $\mathbb{C}$ such that 
$$h_{b}=(G^{-1})^* h^{\beta}.$$

\noindent The Riemannian metric $h_b$ is then periodic with respect to the fundamental domain $[-\frac{\pi}{2},\frac{\pi}{2}]\times[-2b,2b]$. The quotient  $(\mathbb{C},h_b)$ by the subgroup of isometries $<\tau,t_{v'}>$, where $\sigma: z\mapsto \bar{z}+\pi$ and $t_{v'}:z\mapsto z+4ib$, is a Riemannian Klein bottle. Its conformal type is 

$$\beta=\int_{0}^{2b}{\frac{dt}{f(t)}}.\\$$\par

\begin{rem} \label{remsph} If $f$ is equal to the function $cosine$ in a cylinder $C_\alpha=\{(u,v)\in \mathbb{C};\ |v|\leq \alpha\}$ with $\alpha \in ]0,\pi/2[$, then the metric $h_b$ is spherical on $C_\alpha$.
\end{rem} 

\subsection{Singular Riemannian metrics on the Klein bottle}\label{singklein}

We will define singular metrics in the setting of Riemannian polyhedrons. For further details on this notion see \cite{babenko}.\\

A polyhedron is a topological space endowed with a triangulation, i.e.  divided into simplexes glued together by their faces. We denote  by $\Sigma$ an arbitrary simplex of a polyhedron $P$.\\

\begin{defi} 
A Riemannian metric on a polyhedron $P$ is a family of Riemannian metrics $\{g_\Sigma\}_{\Sigma \in I}$, where $I$ is in bijection with the set of simplexes of $P$. 
These metrics should satisfy the following conditions: \begin{enumerate}

\item Every $g_\Sigma$ is a smooth metric in the interior of the simplex $\Sigma$.

\item The metrics $g_\Sigma$ coincide on the  faces; i.e. for any pair of simplexes $\Sigma_1$, $\Sigma_2$, we have the equality
$$g_{\Sigma_1}|_{\Sigma_1 \bigcap \Sigma_2}=g_{\Sigma_2}|_{\Sigma_1 \bigcap \Sigma_2}$$
\end{enumerate}
\end{defi}

\noindent Such a Riemannian structure on the polyhedron  allows us to calculate the length of any piecewise smooth curve in $P$, and makes of the polyhedron $(P,g)$  a length space. The geodesics of a Riemannian polyhedron are the geodesics of the associated length structure (see \cite{bbi}). 
In the interior of a simplex $(\Sigma, g_\Sigma)$, the first variation formula shows that such a geodesic is a geodesic of $g_\Sigma$ in the Riemannian sense (see Prop 2.1 below). Note that the Riemannian measure is defined exactly as in the smooth case.\\

We will introduce now two families of singular Riemannian metrics on the Klein bottle. C. Bavard showed that, for a suitable choice of the parameters, these metrics are systolically extremal in their conformal class (cf. \cite{bavard88}). For every $\mu \in ]0,\infty[$ and every $\lambda \in ]0,\mu[$, define \\

\begin{enumerate}

\item the metric $S_\lambda$ on $\mathbb{C}$ periodic with respect to the fundamental domain $[-\frac{\pi}{2},\frac{\pi}{2}]\times[-2\lambda,2\lambda]$ and defined by 

$$(S_\lambda)_{(u,v)}=s_{\lambda}^{2}(v)du^{2}+dv^{2} $$ where $s_{\lambda}$ is the unique one-variable function invariant by the translation $(u,v) \mapsto (u,v+2\lambda)$ which agrees with  $cosine $ on $[-\lambda,\lambda]$. It is spherical outside the singular lines $v=n \lambda$, where $n \in \mathbb{Z}$. 

\item  the  metric $SF_{\lambda,\mu}$ on $\mathbb{C}$ periodic with respect to the fundamental domain $[-\frac{\pi}{2},\frac{\pi}{2}]\times[-2\lambda,2\lambda]$ and defined by
$$
(SF_{\lambda,\mu})_{(u,v)} =  \sigma_{\lambda,\mu}^{2}(v)du^{2}+dv^{2}$$ 
where $\sigma_{\lambda,\mu}$ is the unique one-variable function invariant by the translation $(u,v) \mapsto (u,v+2\lambda)$ which agrees with  $cosine $ on $[-\mu, \mu]$ and is equal to the constant $\cos(\mu)$ on $[\mu,2\lambda-\mu]$ (see Fig. 1). It is spherical on the band $\mathbb{R} \times [-\mu, \mu]$ and its images by the translations $(u,v)\mapsto (u,v+2n\lambda)$, where $n \in \mathbb{Z}$,  and flat elsewhere. \\
\end{enumerate}
\par
We will denote by $(K,S_\lambda)$ the Riemannian Klein bottle obtained by taking the quotient of $(\mathbb{C},S_\lambda)$ by the subgroup $<\sigma,t_{v'}>$. Its conformal type is $$\beta=2\ln(\tan(\frac{\pi}{4}+\frac{\lambda}{2})).$$
Similarly, $(K,SF_{\lambda,\mu})$ will denote the Riemannian Klein bottle obtained by taking the quotient of $(\mathbb{C},SF_{\lambda,\mu})$ by the subgroup $<\sigma,t_{v'}>$. Its conformal type is
$$\beta=2\ln(\tan(\frac{\pi}{4}+\frac{\mu}{2}))+\frac{1}{\cos(\mu)}(2\lambda-2\mu).$$

\begin{figure}
\begin{center}
\includegraphics[clip, trim=0cm 10cm 9cm 3cm, width=\textwidth]{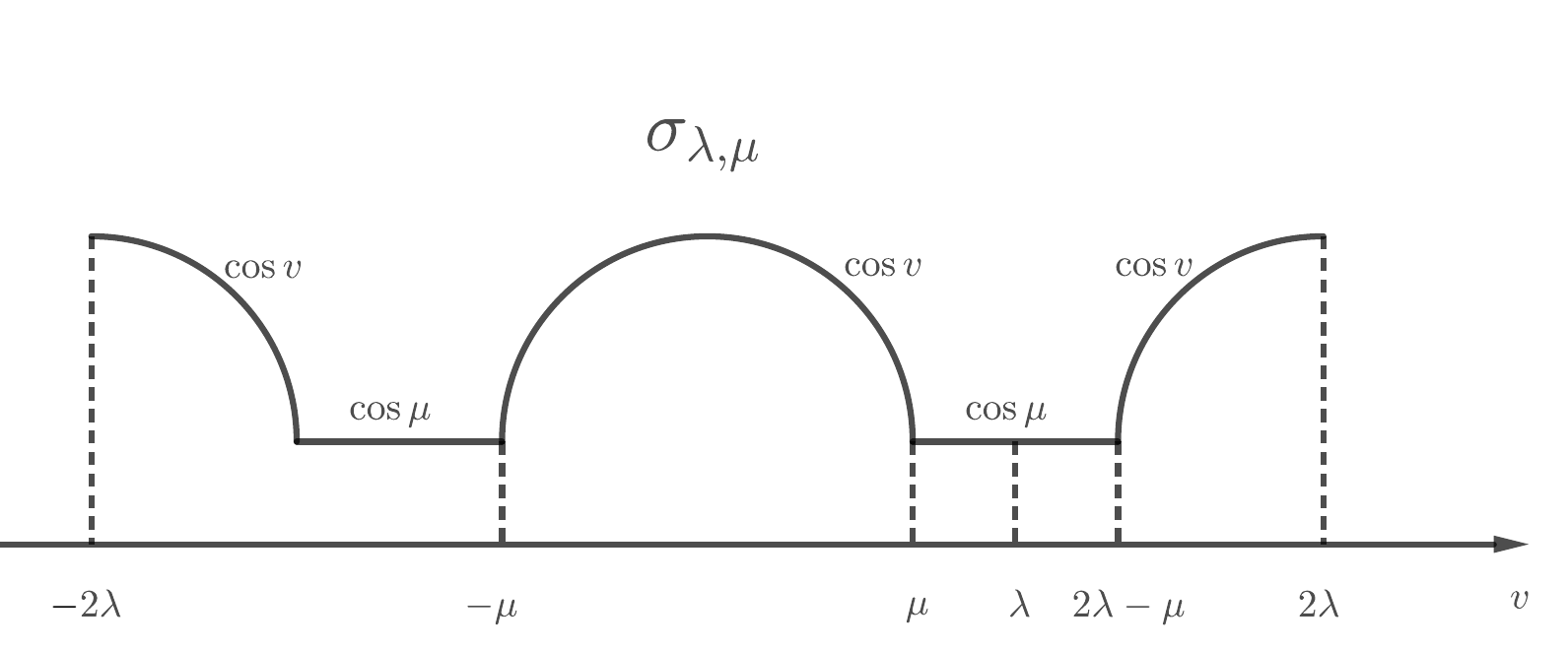}
\caption{The graph of $\sigma_{\lambda,\mu}$}
\end{center}
\end{figure}

In 1986, C. Bavard showed (cf. \cite{bavard86}) that every Riemannian Klein bottle $K$ satisfies the isosystolic inequality

$$   sys(K)^2 \leq \frac{\pi}{2\sqrt{2}} \cdot area(K)$$
where the equality is attained by the singular metric $(K,S_{\frac{\pi}{4}})$ (also denoted $\KK_{\frac{\pi}{4}}$ later in the paper).

\begin{rem}\label{propklein} The isometry group of $(K,S_\lambda)$ and $(K,SF_{\lambda,\mu})$  is of dimension one. Note that an element $r_\delta$ in the connected component of the identity in $Isom(K,SF_{\lambda,\mu})$ corresponds to a rotation of angle $\delta$ in the spherical region of $(K,SF_{\lambda,\mu})$ and a horizontal translation of modulus $\delta\cos \mu$ in the flat region of $(K,SF_{\lambda,\mu})$.
\end{rem}

The geometric properties of $SF_{\lambda,\mu}$ and $S_{\lambda}$ are very similar to those of the extremal Klein bottle $(K,S_{\frac{\pi}{4}})$ (see \cite{elaf08} Section 3). In particular, we have

\begin{prop}\label{geodsing} Let $\gamma$ be a geodesic of $(K,SF_{\lambda,\mu})$ joining a point in the spherical region to a point in the adjacent flat region and let $l$ be the singular line common to both regions. Then $\gamma$ is the juxtaposition of an arc of a great circle in the spherical region and a segment in the flat region such that the sum of the angles between these lines and $l$ is equal to $\pi$.
\end{prop}

\begin{proof}
The proof is very similar to the proof of Proposition 1 in \cite{elaf08}.
\end{proof}

\begin{defi}   For a Riemannian (singular) manifold $(M,g)$, the displacement $d(\gamma)$ of an isometry $\gamma$ of $(M,g)$ is the number

$$d(\gamma)=\inf_{p\in M}dist(p,\gamma(p))$$

\end{defi}


\begin{prop} In $(K,SF_{\lambda,\mu})$, the displacements of $r_\alpha$ and $T_\alpha$ are given by

$$d(r_\delta)= \inf \{\delta \cos \mu,\pi -\delta\}$$

$$d(T_\delta)=\inf\{2\cos^{-1}(\cos (\frac{\delta_1}{2})\cos \mu)+\frac{2(\lambda-\mu)(1-\cos^2 (\frac{\delta_1}{2})\cos^2 \mu)^{\frac{1}{2}}}{\cos (\frac{\delta_1}{2}) \sin \lambda}, \cos \mu(\pi-\delta)\}$$

\noindent where $\delta_1$ is given by the formula $\delta_1+\frac{2(\lambda-\mu)\tan (\frac{\delta_1}{2})}{\sin (\mu) }=\delta$
  
\end{prop}

\begin{proof}
Let $(T,SF_{\lambda,\mu})$ be the covering torus of $(K,SF_{\lambda,\mu})$. Suppose $p$ is a point in the spherical region of $(T,SF_{\lambda,\mu})$ and $p'$ a point in the flat region of $(T,SF_{\lambda,\mu})$. Then, a direct computation of the distance shows that $d_{(T,SF_{\lambda,\mu})}(p,r_\delta(p)) \geq d_{(T,SF_{\lambda,\mu})}(p',r_\delta(p'))=\delta \cos \mu$. On the other hand, the computation of the distance $d_{(T,SF_{\lambda,\mu})}(\tau(p),r_\delta(p))$  shows that it is a decreasing function of the latitude  (see \cite{elaf08} Prop. 3 and Cor. 1 for details). Note that the critical value is $\delta_0= \frac{\pi}{1+\cos \mu}$.\\
Concerning the transformation $T_\delta$, we first note that $d(p,T_\delta(p))$ is constant for all $p \in (T,SF_{\lambda,\mu})$. The proof, left to the reader, uses the same arguments of Prop. 2 in \cite{elaf08}. Now, let $p$ be a point on the singular locus of $(T,SF_{\lambda,\mu})$, we will calculate the distance in $(T,SF_{\lambda,\mu})$ between $p$ and $p'=T_\delta(p)$. By Prop. \ref{geodsing}, the geodesic $\gamma$ joining these two points is the juxtaposition of an arc of a great circle in the spherical region and a segment in the flat region making the same angle $\beta$ with the singular line. We denote by $q$ (resp. $r$) the intersection point of $\gamma$ with the singular line (resp. the equator), and by $p_1$ (resp. $q_1$) the orthogonal projection of $p$ (resp. $q$) on the equator (see Figure 2). If we denote by $\delta_1$ the length of the arc $\widehat{p_1q_1}$, then we have

$$\delta=\delta_1+\frac{2(\lambda-\mu)}{\cos \mu \tan \beta}.$$

Now, in the spherical triangle $prp_1$, the classical formulas of spherical geometry show that

$$\sin \ell (\widehat{pr})=\frac{\sin \frac{\delta_1}{2}}{\cos \beta}\ \textrm{ and }\ \cos \ell (\widehat{pr})= \cos \frac{\delta_1}{2} \cos \mu$$

\noindent where $\widehat{pr}$ denotes the geodesic arc joining $p$ to $r$ in $(T,SF_{\lambda,\mu})$. This shows that $\cos^2 \beta=\frac{\sin^2\frac{\delta_1}{2}}{1-\cos^2 (\frac{\delta_1}{2})\cos^2 \mu}$.
Now, the displacement of $T_\delta$ is equal to the length of $\gamma$, i.e., $d(T_\delta)=2\ell(\widehat{pr})+\ell([qp'])=2\cos^{-1}(\cos (\frac{\delta_1}{2})\cos \mu)+\frac{2(\lambda-\mu)}{\sin \beta}$. Therefore we get  

$$d(T_\delta)=2\cos^{-1}(\cos (\frac{\delta_1}{2})\cos \mu)+\frac{2(\lambda-\mu)(1-\cos^2 (\frac{\delta_1}{2})\cos^2 \mu)^{\frac{1}{2}}}{\cos (\frac{\delta_1}{2}) \sin \lambda}$$

\noindent where $\delta_1+\frac{2(\lambda-\mu)\tan (\frac{\delta_1}{2})}{\sin (\mu) }=\delta$. Finally, the infimum of $d_{(T,SF_{\lambda,\mu})}(\tau(p),T_\delta(p))$ is clearly attained for the points on the flat region of $(T,SF_{\lambda,\mu})$.

\end{proof}


\begin{figure}
\begin{tikzpicture}
    \node[anchor=south west,inner sep=0] (X) at (0,0) {\includegraphics[clip, trim=0cm 5cm 7cm 3cm, width=0.8\textwidth]{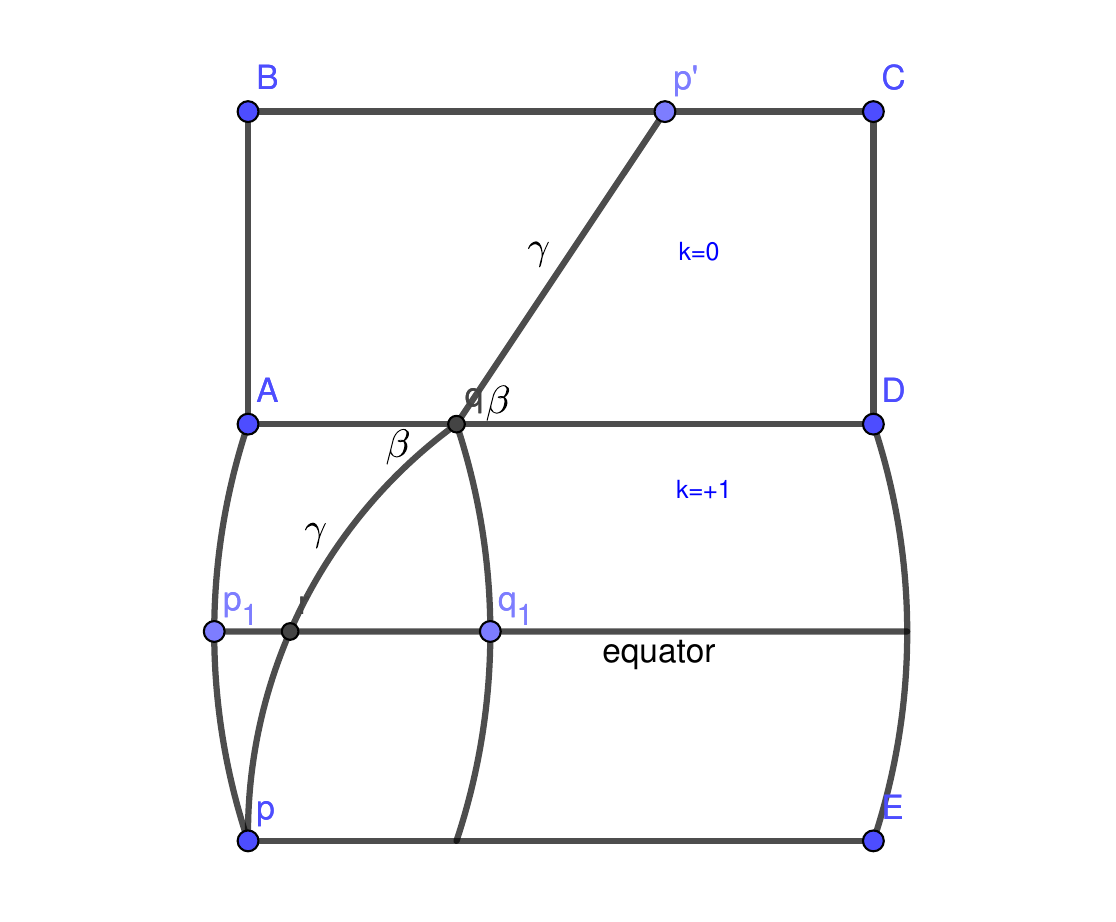}};
    
\node [blue, ultra thick] at ($(X.south west)!0.7!(X.north east)$) {k=0};
\node [blue, ultra thick] at ($(X.east)!0.3!(X.west)$) {k=+1};

\end{tikzpicture}
\caption{The displacement of $T_\delta$ in $(T,SF_{\lambda,\mu})$ }
\end{figure}

\subsection{Bavard's conformal isosystolic inequalities on the Klein bottle}

We can state now Bavard's theorem on the \emph{Conformal Isosystolic Inequality} on the Klein bottle.

\begin{theorem} (\cite{bavard88}) Let $\beta>0$. For every Riemannian metric $g$ on the Klein bottle $K$ of conformal type $\beta$, the following holds

\begin{itemize}

\item If $0 < \beta \leq \frac{\pi}{2}$, then $\frac{Sys^2(g)}{Area(g)} \leq \frac{2\beta}{\pi}$. Moreover, the equality is attained if and only if $g$ is proportional to the flat metric $F_\beta$.

\item If $\frac{\pi}{2} \leq \beta \leq 2 \ln (1+\sqrt{2})$, then $\frac{Sys^2(g)}{Area(g)} \leq \frac{\pi}{4(\sin b+(\frac{\pi}{4}-b)\cos b)}$. Moreover, the equality is attained if and only if $g$ is proportional to the spherical-flat metric $SF_{\frac{\pi}{4},b}$, for $b$ satisfying $2 \ln \tan (\frac{b}{2}+\frac{\pi}{4}) + \frac{\pi}{2}+2b= \beta \cos b$.

\item If $2 \ln (1+\sqrt{2}) \leq \beta \leq 2 \ln (2+\sqrt{3})$, then $\frac{Sys^2(g)}{Area(g)} \leq \frac{\pi}{4\tanh \frac{\beta}{2}}$. Moreover, the equality is attained if and only if $g$ is proportional to the spherical metric $S_{b}$, for $b$ satisfying $2 \ln \tan (\frac{b}{2}+\frac{\pi}{4})= \beta$. 

\item If $2 \ln (2+\sqrt{3}) \leq \beta $, then $\frac{Sys^2(g)}{Area(g)} \leq \frac{2\pi}{4\sqrt{3}+\beta -2\ln (2+\sqrt{3})}$. Moreover, the equality is attained if and only if $g$ is proportional to the spherical-flat metric $SF_{b,\frac{\pi}{3}}$, for $b$ satisfying $2 \ln (2+\sqrt{3})+8(b-\frac{\pi}{3})= \beta$. 

\end{itemize}

\end{theorem}

For each $0< b < +\infty$, we denote

$$ \KK_b=\Big \{ \begin{array}{cc} SF_{\frac{\pi}{4},b} & \textrm{ if } 0 < b \leq \frac{\pi}{4} \\
S_{b} & \textrm{ if } \frac{\pi}{4} \leq b \leq \frac{\pi}{3} \\
SF_{b,\frac{\pi}{3}} & \textrm{ if } \frac{\pi}{3} \leq b < +\infty \end{array}$$

\begin{rem} Note the shifting role of the parameter $b$: it corresponds to the maximum latitude of the spherical band if $0\leq b \leq \frac{\pi}{3}$ and to half the width of the flat band (minus $\frac{\pi}{3}$)  if $\frac{\pi}{3} \leq b < +\infty$.
\end{rem}

\section{Flat non-orientable $3$-manifolds and the singular metric $g_i^{b,c}$}\label{g0}

 The classification of flat manifolds of dimension $3$ results of a direct method of classification of discrete, cocompact subgroups of $\mathrm{Isom}(\RR^3)$ operating freely. 
This classification is due to W. Hantzsche and H. Wendt  (1935), and exposed in the book \cite{wolf} of J.A. Wolf.
There exist ten compact and flat manifolds of dimension $3$ up to an affine diffeomorphism. In the non-orientable case, there are four types of compact, flat $3$-manifolds (see \cite{wolf} and \cite{elaf08}). They can be obtained as a twisted product of a flat Klein bottle by an isometry of $\mathrm{Isom}(\KK)/\mathrm{Isom}_0(\KK)=\ZZ_2 \times \ZZ_2$ (cf. \cite{wolf}, \cite{elaf08}).

We define a singular two-parameter family of Riemannian metrics $g_i^{b,c}$ on each one of the four types of non-orientable Bieberbach manifolds. These metrics have ``a lot" of systolic geodesics. For each $0\leq b < +\infty$, we denote

\begin{enumerate}

\item $(B_1,g_1^{b,c})$: it is the quotient of $\KK_b \times \RR$ by the subgroup generated by $(p,t)\mapsto \big(r_\delta(p), t+c\big)$ where $\delta$ is the unique real number that maximizes the the displacement of $r_{\delta}$ in $\KK_b$ (see Remark \ref{remklein} below) and  $c\geq \sqrt{\pi^2-\big(\mathrm{d}(r_{\delta})\big)^2}$. The range of values of $c$ guarantees that  the geodesics closed by the isometry $(r_\delta,t_c)$ are of length greater than $\pi$ (hence $sys(B_1,g_1^{b,c})=\pi$).\\

\item $(B_2,g_2^{b,c})$: it is the quotient of $\KK_b\times \RR$ by the subgroup generated by $(p,t)\mapsto \Big(T_{\delta}(p), t+c\Big)$ where $\delta$ is the unique real number that maximizes the displacement of $T_{\delta}$ in $\KK_b$ (see Remark \ref{remklein} below) and  $c\geq \sqrt{\pi^2-\big(\mathrm{d}(T_{\delta})\big)^2}$. The range of values of $c$ guarantees that  the geodesics closed by the isometry $(T_\delta,t_c)$ are of length greater than $\pi$ (hence $sys(B_2,g_2^{b,c})=\pi$).\\

\item $(B_3,g_3^{b,c})$: it is the quotient of $\KK_b \times \RR$ by the subgroup generated by $(p,t))\mapsto \big(S_1(p), t+c \big)$, where $c\geq \pi$. Note that the geodesics closed by the isometry $(S_1,t_c)$ are of length greater than $\pi$ (hence $sys(B_3,g_3^{b,c})=\pi$).\\

\item $(B_4,g_4^{b,c})$: it is the quotient of $\KK_b \times \RR$ by the subgroup generated by $(p,t))\mapsto \big(S_2(p), t+c \big)$. Note that the geodesics closed by the isometry $(S_2,t_c)$ are of length greater than $\pi$ (hence $sys(B_4,g_4^{b,c})=\pi$).\\

\end{enumerate}

\begin{rem} \label{remklein} In the construction of the manifold $(B_1,g_1^{b,c})$, the value of $\delta$ which maximizes the displacement of  $r_{\delta}$ is

$$\delta=\left \{\begin{array}{ccc}\frac{\pi}{(1+\cos b)} & if & 0<b \leq \frac{\pi}{3}\\
\frac{2\pi\sqrt{2}}{3} & if & b > \frac{\pi}{3}\end{array}\right.$$ 


In the construction of the manifold $(B_2,g_2^{b,c})$, the value of $\delta$ which maximizes the displacement of  $T_{\delta}$ is (recall the formula $\delta_1+\frac{2(\lambda-\mu)\tan (\frac{\delta_1}{2})}{\sin (\mu) }=\delta$ from Prop. \ref{geodsing})

$$\delta=\left \{\begin{array}{ccc}2\cos^{-1}(\cos (\frac{\delta_1}{2})\cos b)+\frac{2(\frac{\pi}{4}-b)(1-\cos^2 (\frac{\delta_1}{2})\cos^2 b)^{\frac{1}{2}}}{\cos (\frac{\delta_1}{2}) \frac{\sqrt{2}}{2}}=\cos b (\pi-\delta) & if & 0<b \leq \frac{\pi}{4}\\
2\cos^{-1}(\cos (\frac{\delta}{2})\cos b)= \cos b (\pi-\delta) & if & \frac{\pi}{4}<b \leq \frac{\pi}{3}\\
2\cos^{-1}(\cos (\frac{\delta_1}{2})/2)+\frac{2(b-\frac{\pi}{3})(1-\cos^2 (\frac{\delta_1}{2})/4)^{\frac{1}{2}}}{\cos (\frac{\delta_1}{2}) \sin b}=(\pi-\delta)/2  & if &\frac{\pi}{3}<b < +\infty
\end{array}\right. $$

\end{rem} 

\section{Proof of theorem 1}\label{proof1}

First we note that the metrics $g_i^{b,c}$ can be written in the form

$$g_i^{b,c}=\Big \{ \begin{array}{cc} \sigma_{\frac{\pi}{4},b}^{2}(v)du^{2}+dv^{2}+dw^2 & \textrm{ if } 0\leq b \leq \frac{\pi}{4} \\
s_{b}^{2}(v)du^{2}+dv^{2}+dw^2 & \textrm{ if } \frac{\pi}{4} \leq b \leq \frac{\pi}{3} \\
\sigma_{b,\frac{\pi}{3}}(v)du^{2}+dv^{2}+dw^2 & \textrm{ if } \frac{\pi}{3} \leq b < +\infty \end{array}$$ 

\noindent where $s$ and $\sigma$ are defined in Section \ref{singklein}.\\

In the following, we will prove Theorem 1. By Theorem \ref{confext}, it is sufficient to construct a measure $\mu$ on a certain set $S$ of systolic geodesics of $(B_i,g_i^{b,c})$ ($1\leq i \leq 4$) such that:\\

$$^*\mu=\frac{Sys(g_0)}{Vol(g_0)} \cdot dg_i^{b,c}$$

where $dg_i^{b,c}$ is the volume measure of $(B_i,g_i^{b,c})$.\\

We define first the family $F_1(b)$ (where $0<b< +\infty$) of systolic geodesics as follows.\\

\noindent Consider for each $\theta \in \mathbb{R}/\pi \mathbb{Z}$ and each $a \in [-b,b] \cup [\frac{\pi}{2}-b, \frac{\pi}{2}]\cup [-\frac{\pi}{2},-\frac{\pi}{2}+b]$, the geodesic loop $\alpha_{\theta,a,w}$ (image in the hypersurface $\{w=w_0\}$ in $(B_i,g_i^{b,c})$ of a great circle) going through the points $(\theta-\pi/2,0,w)$ and $(\theta,a,w)$. Then we set  \\

$$F_{1}(b)= \left\{\alpha_{\theta,a,w} \mid  |a|\leq b \textrm{ or } \frac{\pi}{2}-b \leq |a|\leq \frac{\pi}{2}, \theta\in \mathbb{R}/\pi \mathbb{Z}\right\}.$$

\noindent The proof of C. Bavard in \cite{bavard88} works prefectly in the three dimensional case for each one of the four manifolds $(B_i,g_i^{b,c})$.\\

\begin{enumerate}

\item Let $0\leq b \leq \frac{\pi}{4}$.\\

\noindent Define the family $F_2$ of systolic geodesics as follows\\

$$F_2=\{\gamma_{u_0,w_0}:\ \gamma_{u_0,w_0}(t)=(u_0,t,w_0)\}$$ where $u_0 \in \RR/\pi\ZZ$ and $0\leq w_0 \leq c$.\\

\noindent We consider in each one of the manifolds $(B_i,g_i^{b,c})$ ($i\in \{1,2,3,4\}$), both families of geodesics $F_1(b)$ and $F_2$.

\noindent We equip the family $F_{2}$ with the measure 
$$\mu_{2} = \cos(b)du\otimes dw.$$

\noindent Then, it can be easily verified that 
$$
^*\mu_2 =  \frac{\cos(b)}{\sigma_{\frac{\pi}{4},b}(v)} \cdot d(g_i^{b,c}). $$ 

\noindent Next, we define\\

$$ \begin{array}{cccc} h: & [-b,b] & \longrightarrow & \mathbb{R}\\
& a & \mapsto &\frac{\sin(|a|)}{\pi\cos(a)}\sqrt{\cos^{2}(a)-\cos^{2}(b)}. \\ \end{array}$$ 

\noindent Now, we equip the family $F_{1}(b)$ with the measure
\begin{equation*}
\mu_{1} ={h}(a)\ da\otimes d\theta \otimes dw
\end{equation*}

\noindent  Moreover, an adaptation of Lemma 2.5 in \cite{elya} to the three dimensional case gives

$$\chi_{\{|v|\leq b\}}\ ^*\mu_2=2\chi_{\{|v|\leq b\}}\Big (\int^{b}_{|v|}\left(\cos^{2}(v)-\cos^{2}(a)\right)^{-\frac{1}{2}}{h}(a)\ da\Big ) \cdot d(g_i^{b,c})	$$

\noindent Then, we derive 

$$\chi_{\{|v|\leq b\}}\ ^*\mu_2 =\chi_{\{|v|\leq b\}}\  \big(1-\frac{\cos(b)}{\sigma_{\frac{\pi}{4},b}(v)}\big) \cdot d(g_i^{b,c}).$$ 

\noindent Finally, since $\sigma_{\frac{\pi}{2},b}$ and ${h}$ are invariant by the translation $v \mapsto v+\frac{\pi}{2}$, we have

$$^*\mu_2 =\big(1-\frac{\cos(b)}{\sigma_{\frac{\pi}{4},b}(v)}\big) \cdot d(g_i^{b,c}).$$

\noindent Hence, we get $^*\mu_1+ ^*\mu_2=d (g_i^{b,c})$. 
Finally, from Theorem \ref{confext}, we deduce that  $g_i^{b,c}$ is extremal in its conformal class.\\

\item Let $\frac{\pi}{4} \leq b \leq \frac{\pi}{3}$. We consider in each one of the manifolds $(B_i,g_i^{b,c})$ ($i\in \{1,2,3,4\}$), the family of geodesics $F_1(b)$.\\

\noindent We define\\

$$ \begin{array}{cccc} h: & [-b,b] & \longrightarrow & \mathbb{R}\\
& a & \mapsto &\frac{\sin(|2a|)}{2\pi}({\cos^{2}(a)-\cos^{2}(b)})^{-\frac{1}{2}} \\ \end{array}$$ 

\noindent We equip the family $F_{1}(b)$ with the measure
\begin{equation*}
\mu ={h}(a)\ da\otimes d\theta \otimes dw
\end{equation*}

\noindent Therefore, we get

$$\chi_{\{|v|\leq b\}}\ ^*\mu=2\chi_{\{|v|\leq b\}}\Big (\int^{b}_{|v|}\left(\cos^{2}(v)-\cos^{2}(a)\right)^{-\frac{1}{2}}{h}(a)\ da\Big ) \cdot d(g_i^{b,c})	$$

\noindent Then, we derive 

$$\chi_{\{|v|\leq b\}}\ ^*\mu =\chi_{\{|v|\leq b\}}\  \cdot d(g_i^{b,c}).$$ 

\noindent Since $s_{b}$ and ${h}$ are invariant by the translation $v\mapsto v+2b$,  we get $^*\mu=d (g_i^{b,c})$. 
Finally, from Theorem \ref{confext}, we deduce that  $g_i^{b,c}$ is extremal in its conformal class.\\

\item Let $\frac{\pi}{3} \leq b < +\infty$. \\

\noindent Define the family $F_2$ of systolic geodesics as follows\\

$$F_2=\{\gamma_{v_0,w_0}:\ \gamma_{v_0,w_0}(t)=(t,v_0,w_0)\}$$ 

\noindent where $v_0\in [-2b+\frac{\pi}{3},-\frac{\pi}{3}] \cup [\frac{\pi}{3},2b-\frac{\pi}{3}]$ and  $0\leq w_0 \leq c$.\\

\noindent We consider in each one of the manifolds $(B_i,g_i^{b,c})$ ($i\in \{1,2,3,4\}$), both families of geodesics $F_1(\frac{\pi}{3})$ and $F_2$.

\noindent We equip the family $F_{2}$ with the measure 
$$\mu_{2} = dv\otimes dw.$$

\noindent Then, it can be easily verified that 
$$
^*\mu_2 =  \frac{{1}}{\sigma_{b,\frac{\pi}{3}}(v)} \cdot d(g_i^{b,c}) $$ 

\noindent Next, we define\\

$$ \begin{array}{cccc} h: & [-\frac{\pi}{3},\frac{\pi}{3}] & \longrightarrow & \mathbb{R}\\
& a & \mapsto &\frac{\sin(|2a|)}{2\pi}(\cos^{2}(a)-\frac{1}{2})^{-\frac{1}{2}} \\ \end{array}$$ 

\noindent Now, we equip the family $F_{1}(\frac{\pi}{3})$ with the measure
\begin{equation*}
\mu_{2} ={h}(a)\ da\otimes d\theta \otimes dw
\end{equation*}

\noindent Therefore, we get

$$\chi_{\{|v|\leq b\}}\ ^*\mu_1=2\chi_{\{|v|\leq b\}}\Big (\int^{b}_{|v|}\left(\cos^{2}(v)-\cos^{2}(a)\right)^{-\frac{1}{2}}{h}(a)\ da\Big ) \cdot d(g_i^{b,c})	$$

\noindent Then, we derive 

$$\chi_{\{|v|\leq b\}}\ ^*\mu_1 =\chi_{\{|v|\leq b\}}\  \big(1-\frac{1}{\sigma_{b,\frac{\pi}{3}}(v)}\big) \cdot d(g_i^{b,c}).$$ 

\noindent Finally, since $\sigma_{b,\frac{\pi}{3}}$ and ${h}$ are invariant by the translation $(u,v)\mapsto (u,v+2b)$, we have

$$^*\mu_1 =\big(1-\frac{1}{\sigma_{b,\frac{\pi}{3}}(v)}\big) \cdot d(g_i^{b,c}).$$

\noindent Hence, we get $^*\mu_1+ ^*\mu_2=d (g_i^{b,c})$. 
Finally, from Theorem \ref{confext}, we deduce that  $g_i^{b,c}$ is extremal in its conformal class.

\end{enumerate}

\section{Maximizing the systolic ratio on the set of parameters}\label{maximus}

In this final section, we calculate, for each $1\leq i \leq 4,$ the maximum of the systolic ratio  $SR_i(b,c)=\frac{sys^3(g_i^{b,c})}{vol(g_i^{b,c})}$  of $(B_i,g_i^{b,c})$ over the set of the parameters $b$ and $c$.\\

Since for each $1\leq  i\leq 4$, the systole of $(B_i,g_i^{b,c})$ is equal to $\pi$ we have

$$SR_i(b,c)=\frac{\pi^3}{Area(\KK_b)\cdot c}$$

Then, we get $$SR_i(b,c)=\left\{ \begin{array}{ccc} \frac{\pi^2}{4(\sin b+\cos (\frac{\pi}{4}-b))\cdot c} & if & 0<b\leq \frac{\pi}{4} \\
 \frac{\pi^2}{4 \sin b\cdot c} & if & \frac{\pi}{4}\leq b \leq \frac{\pi}{3} \\
\frac{\pi^2}{2(\sqrt{3}+(b-\frac{\pi}{3}))\cdot c} & if & \frac{\pi}{3}\leq b
\end{array}\right.$$

\noindent where $c$ is bounded from below by $\pi$ if $i=3,4$, by $\sqrt{\pi^2-\big(\mathrm{d}(r_{\delta})\big)^2}$ if $i=1$ and by $\sqrt{\pi^2-\big(\mathrm{d}(T_{\delta})\big)^2}$ if $i=2$ (see Section \ref{g0} for the expression of the lower bound). It is easy to find the maximum of $SR_3(b,c)$ and $SR_4(b,c)$ since the bound of $c$ does not depend on $b$. Therefore, to maximize $Area(\KK_b)\cdot c$, one should take $c=\pi$ and $b=\frac{\pi}{4}$. This shows that the maximum of $SR_3(b,c)$ and $SR_4(b,c)$ are attained precisely by the metrics obtained as a twisted product of the globally extremal Klein bottle $\KK_{\frac{\pi}{4}}$. The cases of $SR_1(b,c)$ and $SR_2(b,c)$ are a bit more involved. The maximum of these systolic ratios can be calculated using a numerical calculations software (MATLAB for example). It turns out that the maximum is also attained for the metrics obtained as a twisted product of $\KK_{\frac{\pi}{4}}$. More precisely, $\sup SR_1(b,c)$ is attained for $b=\frac{\pi}{4}$ and $c=\pi(2\sqrt 2 -2)^\frac{1}{2}$ and $\sup SR_2(b,c)$ is attained for $b=\frac{\pi}{4}$ and $c= \sqrt{\pi^2-\frac{(\pi-\delta_0)^2}{2}}$, where $\delta_0$ is the solution of the equation $\frac{\cos \delta_0 - 1}{2}=\cos \frac{\pi-\delta_0}{\sqrt{2}}$. An explicit calculation of $\sup SR_i(b,c)$ shows that  the corresponding metrics $(B_i,g_i^{b,c})$ (obtained for the optimal parameters $b$ and $c$) beat every flat Riemannian metric on $B_i$ as it was shown in \cite{elaf08}. Moreover, these metrics satisfy some``nice" geometric properties

\begin{itemize}

\item The systolic geodesics cover the manifold.

\item For any point outside a set of null measure, there exists an infinite number of systolic geodesics going through the point.

\end{itemize}

We summarize the previous discussion with the following table

\begin{table}[htbp]\renewcommand{\arraystretch}{1.2}
\begin{center}
\begin{tabular}[b]{|*{3}{c|}}
   \hline
     Type    &      $\sup{SR_i(b,c)}$                                        &       corresp. values of $b$ and $c$\\
   \hline
   $B_1$    &    $\frac{\pi}{4(\sqrt{2}-1)^{\frac{1}{2}}}$    &  $b=\frac{\pi}{4}$ and $c=\pi(2\sqrt 2 -2)^\frac{1}{2}$\\
   \hline
   $B_2$   &       $\frac{\pi2}{2d_0\sqrt{2}}$                       &    $b=\frac{\pi}{4}$ and $c= \sqrt{\pi^2-\frac{(\pi-\delta_0)^2}{2}}$\\
   \hline
   $B_3$   &      $\frac{\pi}{2\sqrt{2}}$                               &      $b=\frac{\pi}{4}$ and $c=\pi$  \\
   \hline
   $B_4$    &    $\frac{\pi}{2\sqrt{2}}$                               &       $b=\frac{\pi}{4}$ and  $c=\pi$ \\
   \hline
   \end{tabular}\end{center}
\end{table}

Finally, we mention that our result can be generalized to Bieberbach manifolds that can be written as a twisted product of a Klein bottle by an isometry. More precisely, if $B$ is a Bieberbach n-manifold which is diffeomorphic to the quotient of $\KK\times \RR^{n-2}$ (equipped with a flat Riemannian metric) by an isometry of the form $(f,t_a)$, where $f\in \mathrm{Isom}(\KK)/\mathrm{Isom}_0(\KK)$ and $t_a $ is a translation of vector $a\in \RR^{n-2}$,  then $B$ admits a {n-1}-parameter family of singular (non constant curvature) Riemannian metrics each of which is extremal in its conformal class. We believe that the proof works exactly the same way as in the three dimensional case.

\sc{Departement de math\'ematiques, Universit\'e Libanaise, Tripoli, Liban.}


\end{document}